\documentclass{amsart}

\usepackage{amsmath,amssymb,amsthm}
\usepackage{graphicx,tikz}

\newtheorem*{thm}{Theorem}

\theoremstyle{definition}

\theoremstyle{remark}

\begin{document}

\title[]{Using Expander Graphs to test \\whether samples are i.i.d.}
\subjclass[2010]{} 
\thanks{S.S. is supported by the NSF (DMS-1763179) and the Alfred P. Sloan Foundation.}

\author[]{Stefan Steinerberger}
\address{Department of Mathematics, University of Washington, Seattle}
\email{steinerb@uw.edu}

\begin{abstract} The purpose of this note is to point out that the theory of expander graphs leads to an interesting test whether $n$ real numbers $x_1, \dots, x_n$ could be $n$ independent samples of a random variable. To any distinct, real numbers $x_1, \dots, x_n$, we associate a 4-regular graph $G$ as follows: using $\pi$ to denote the permutation ordering the elements, $x_{\pi(1)} < x_{\pi(2)} < \dots < x_{\pi(n)}$, we build a graph on $\left\{1, \dots, n\right\}$ by connecting $i$ and $i+1$ (cyclically) and $\pi(i)$ and $\pi(i+1)$ (cyclically). If the numbers are i.i.d. samples, then a result of Friedman implies that $G$ is close to Ramanujan. This suggests a test for whether these numbers are i.i.d: compute the second largest (in absolute value) eigenvalue of the adjacency matrix.  The larger $\lambda - 2\sqrt{3}$, the less likely it is for the numbers to be i.i.d. We explain why this is a reasonable test and give many examples. 
\end{abstract}
\maketitle

\vspace{-12pt}

\section{Introduction}
\subsection{Introduction.} Let $x_1, \dots, x_n$ be a set of $n$ distinct, real numbers. We will associate to them a unique 4-regular graph as follows: we connect $x_{i}$ to $x_{i+1}$ (cyclically, so $x_n$ also gets connected to $x_1$). Then we order them in increasing size $x_{\pi(1)} < x_{\pi(2)} < \dots < x_{\pi(n)}$ and connect $x_{\pi(i)}$ to $x_{\pi(i+1)}$ and conclude by also connecting to $x_{\pi(n)}$ to $x_{\pi(1)}$. This results in a 4-regular graph.


\begin{center}
\begin{figure}[h!]
\begin{tikzpicture}
\node at (0,0) {\includegraphics[width=0.5\textwidth]{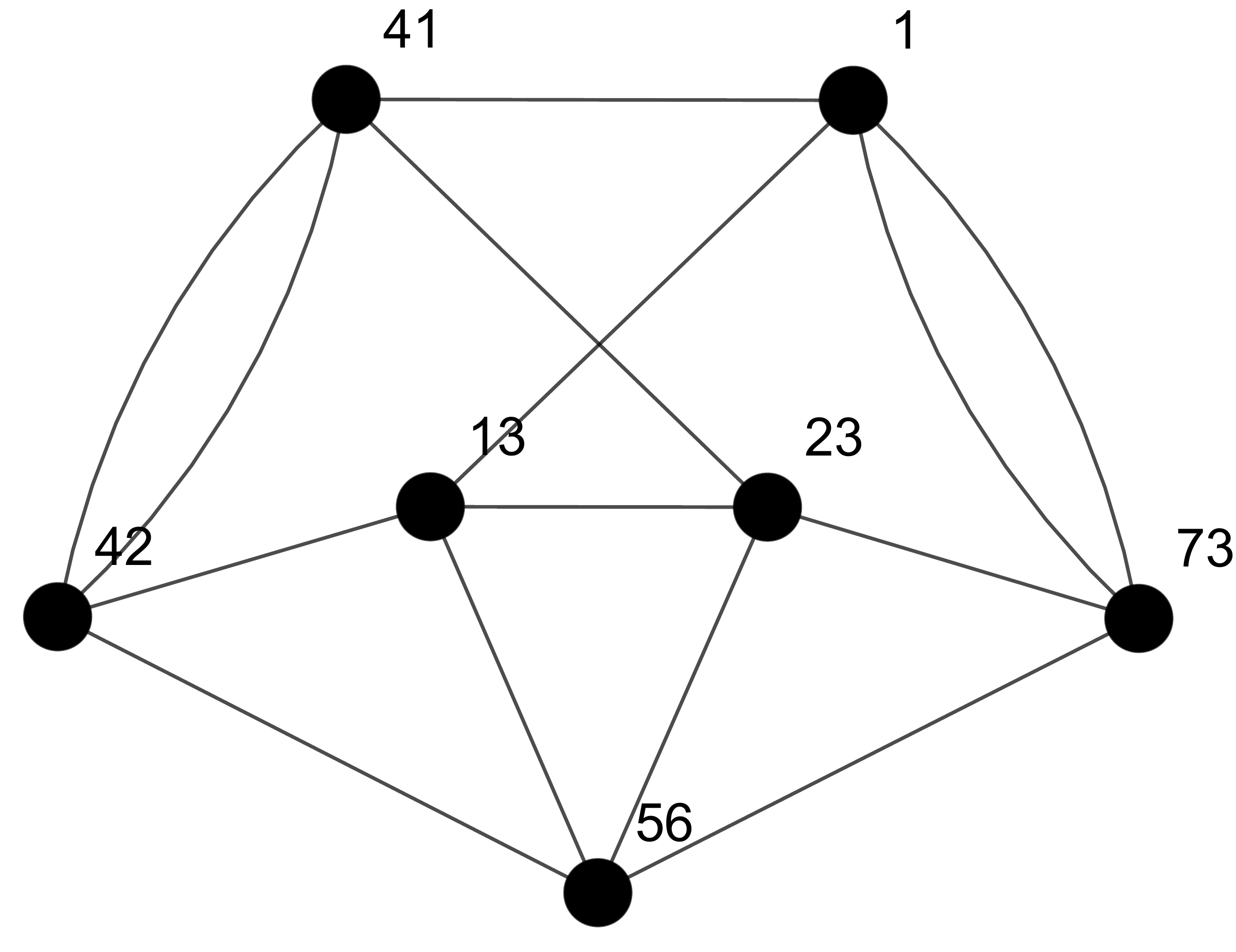}};
\end{tikzpicture}
\caption{The Graph for the first few digits of $\sqrt{2}$: $1, 41, 42, 13, 56, 23, 73$. If the digits of $\sqrt{2}$ behave like truly random numbers, these graphs have optimal expansion properties.}
\end{figure}
\end{center}

We are interested in properties of this type of graph when $x_1, \dots, x_n$ are independent random variables drawn from the same distribution. We assume, for simplicity, that the distribution does not have any atomic mass, the likelihood of having two reals assume the same value is 0. 
This type of graph may be quite interesting for deterministic sequences that are far from random. In fact, we do observe some quite interesting patterns for various deterministic sequences. Fig. 2 shows the type of graph obtained from the van der Corput sequence in base 2
$$ \frac{1}{2}, \frac{1}{4}, \frac{3}{4}, \frac{1}{8}, \frac{5}{8}, \frac{3}{8}, \frac{7}{8}, \frac{1}{16}, \dots$$
and the van der Corput in base 3. In both cases, they seem to approximate an underlying topological limiting object. Something similar seems to happen when building such graphs using irrational multiples such as $x_n = n \sqrt{2} ~\mod 1$. 

\begin{center}
\begin{figure}[h!]
\begin{tikzpicture}
\node at (0,0) {\includegraphics[width=0.3\textwidth]{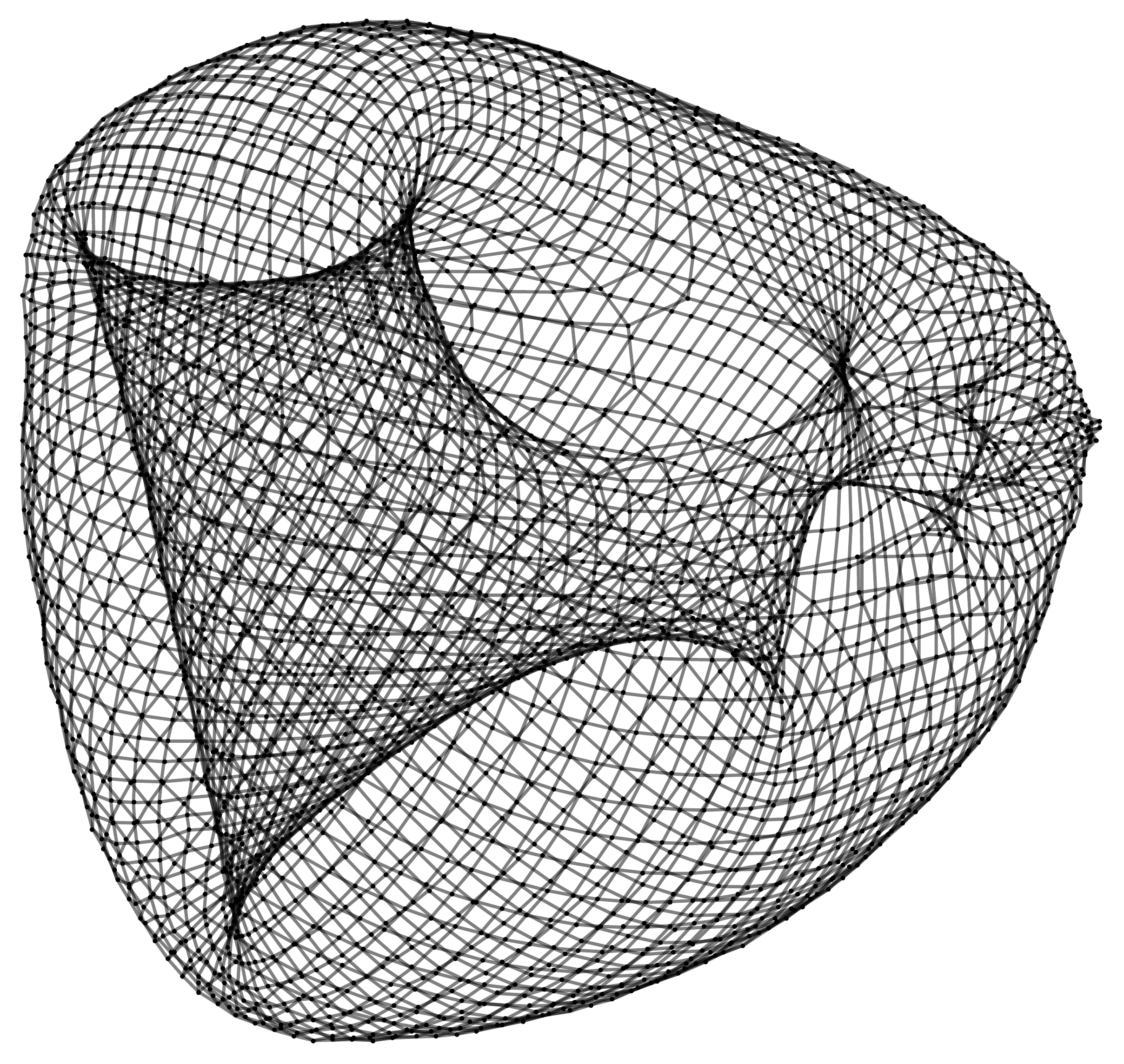}};
\node at (6,0) {\includegraphics[width=0.3\textwidth]{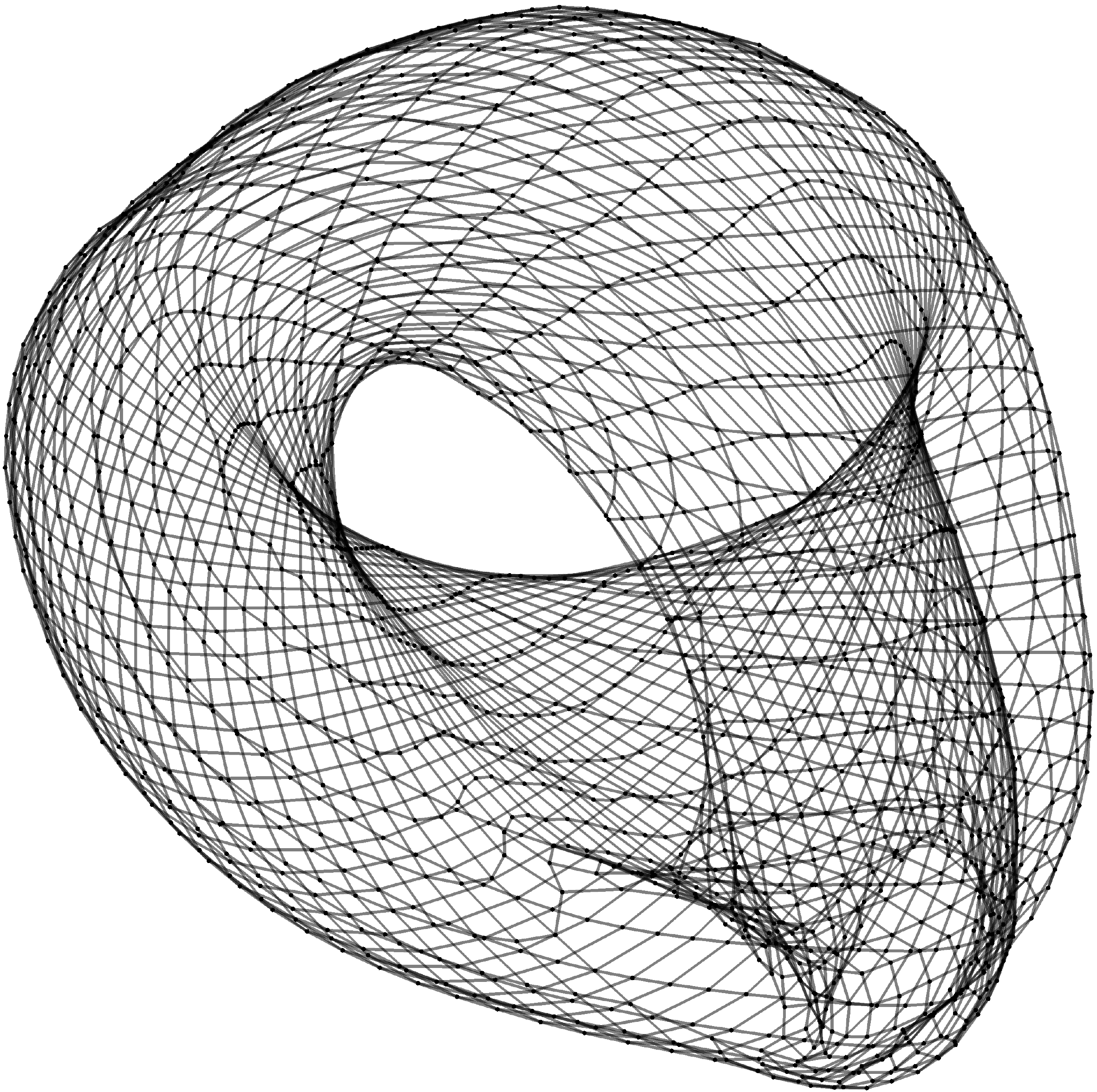}};
\end{tikzpicture}
\caption{The first 4096 elements of the van der Corput sequence in base 2 (left) and the first $2187=3^7$ elements of the van der Corput sequence in base 3 (right).}
\end{figure}
\end{center}

A natural question is whether this graph can be used to understand whether a sequence is random or not. We illustrate this with two more examples: the first is the sequence of Fibonacci number $F_n$ considered in the residue class 10001. The second example is derived from using 100 `random' real numbers provided by a computer.
The picture is different: the first Fibonacci numbers are below 10001, so $x_n$ is strictly increasing for a while (this is the isolated region in the Graph). 

\begin{center}
\begin{figure}[h!]
\begin{tikzpicture}
\node at (0,0) {\includegraphics[width=0.4\textwidth]{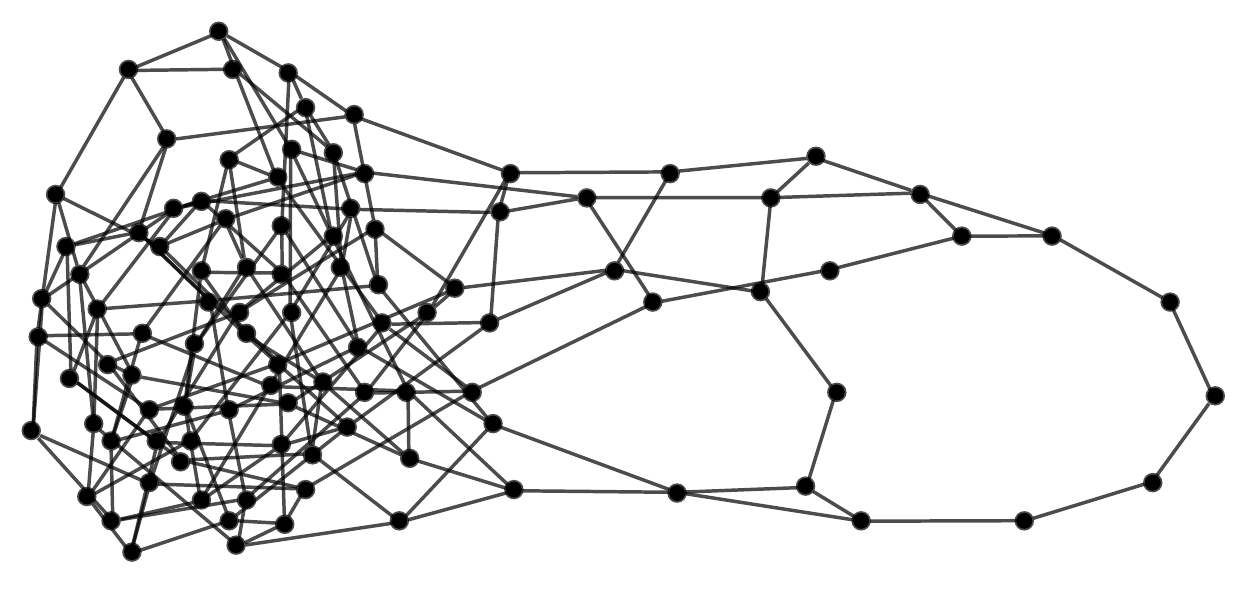}};
\node at (6,0) {\includegraphics[width=0.3\textwidth]{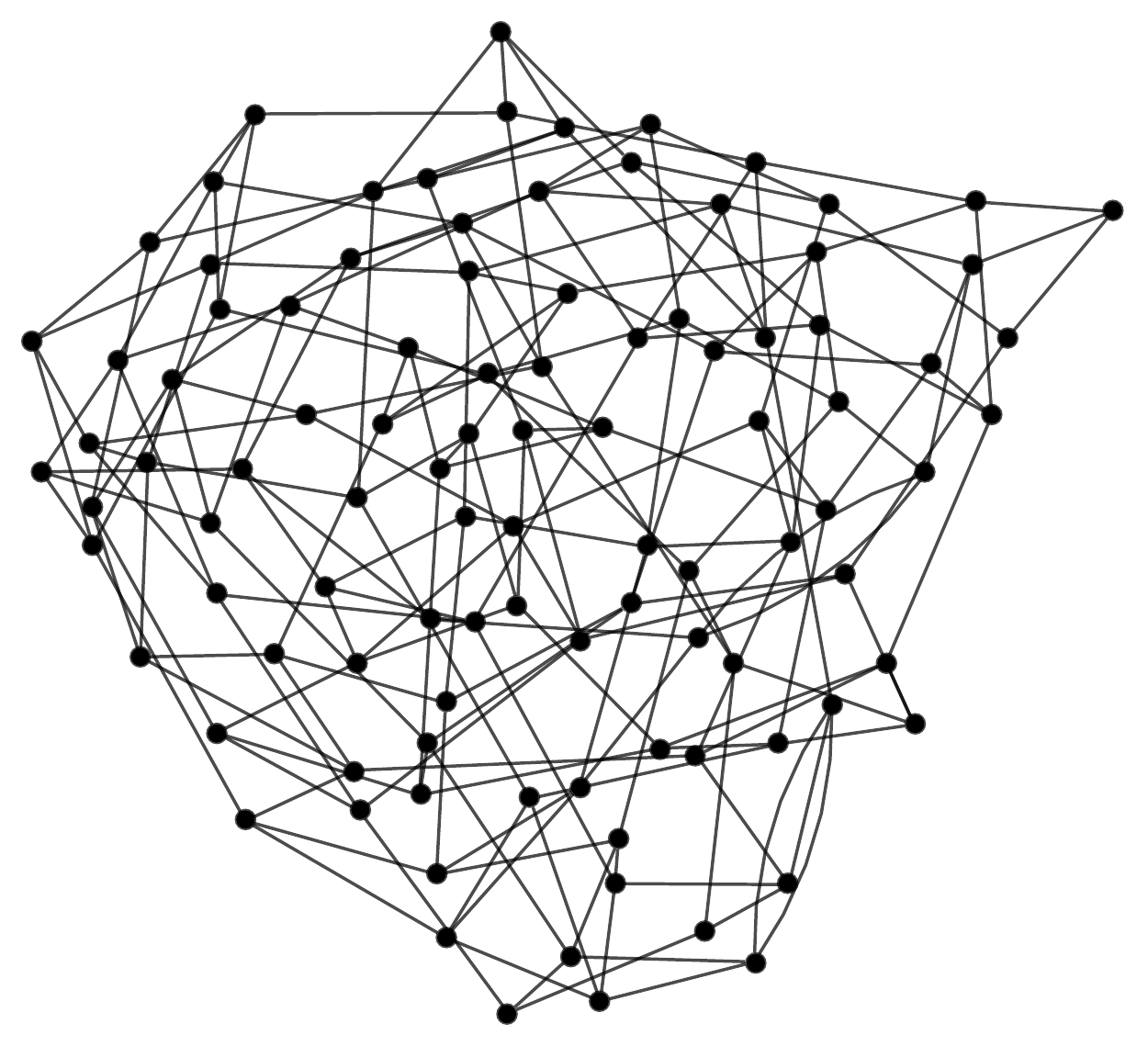}};
\end{tikzpicture}
\caption{$x_n = (F_n\mod 10001)_{n=2}^{101}$ (left) and 100 random reals (right).}
\end{figure}
\end{center}

\subsection{Expanders.}  Expanders are graphs with the property that there are no `insulated' subsets: every subset of certain sizes has a lot of edges connecting it to its complement. Formally, we define the \textbf{expansion ratio} of a graph $G=(V,E)$ as
$$ h(G) = \min_{S \subset V\atop |S| \leq |V|/2} \frac{|\partial S|}{|S|},$$
where $\partial S$  is the set of all edges between $S$ and $V \setminus S$. A sequence $(G_i)_{i=1}^{\infty}$  of $d-$regular graphs with increasing size increasing is a family of expander graphs if the expansion ratio is uniformly bounded away from 0. There are several mostly equivalent points of view: expander graphs are graphs that have the property that despite the number of edges $|E|$ being rather small (comparable to the number of vertices $|V|$), the random walk equidistributes very quickly. 

\begin{center}
\begin{figure}[h!]
\begin{tikzpicture}
\filldraw (0,0) circle (0.05cm);
\filldraw (1,1) circle (0.05cm);
\filldraw (0.2,0.8) circle (0.05cm);
\filldraw (0.7,0) circle (0.05cm);
\draw [thick] (0,0) -- (1,1) -- (0.7,0) -- (0.2, 0.8) -- (0,0);
\draw[thick] (0.5, 0.5) circle (1cm);
\draw [thick, dashed] (0,0) -- (-0.9, -0.2);
\draw [thick, dashed] (0,0) -- (0.9, -0.6);
\draw [thick, dashed] (0.2,0.8) -- (-1, 1);
\draw [thick, dashed] (1,1) -- (2.3, 0.8);
\draw [thick, dashed] (0.7,0) -- (2, 0.2);
\end{tikzpicture}
\caption{Each subset $S$ has a number of edges to the complement that is proportional to $|S|$ (up to $|S| \leq |V|/2$ to ensure there is reasonably sized complement).}
\end{figure}
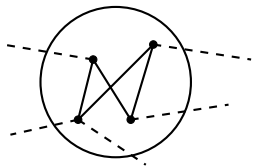
\end{center}

The first proof of the existence of expanders is due to Pinkser \cite{pinsker} in 1973 (however, see also Kolmogorov \& Barzdin \cite{kol}). A typical random $d-$regular graph is known to be an expander but it is actually fairly difficult to give explicit description of expanders. Several constructions are now known
\cite{bo, bgt3, k2, lps1, lsv2, m1, m2, m3, mo, ps2, rvw, va2} -- for this and other aspects, we refer to excellent expository material provided by Hoory, Linial \& Widgerson \cite{linial}, Lubotzky \cite{lub} and Kowalski \cite{kow}.

\subsection{The Spectral Point of View.} There's a spectral point of view on expansion that will be very useful. Given a simple $d-$regular graph $G$, we can consider its adjacency matrix $A \in \mathbb{R}^{|V| \times |V|}$ where
$$ A_{ij} = \begin{cases} 1 \qquad &\mbox{if}~i \sim_E j \\
0 \qquad &\mbox{otherwise.} \end{cases}$$
If $G$ has multipled edges, we replace 1 by the number of edges. This symmetric graph has $|V|$ eigenvalues of which the largest is $d$ (the associated eigenvector is the constant vector). Ordering the eigenvalues $d= \lambda_1 > |\lambda_2| \geq \dots \geq |\lambda_n|$, we see that a relevant quantity is how close such an eigenvalue (different from $\lambda_1$) can get to $d$: the larger the gap, the quicker the random walk converges to the uniform distribution. We first observe that there is a classical connection between the isoperimetric ratio $h(G)$ and the gap (see e.g. \cite{alo86, am85, bus, che, dod84})
$$ \frac{d- |\lambda_2|}{2} \leq h(G) \leq \sqrt{2d(d-|\lambda_2|)}.$$
This inequality thus leads to an alternative and purely spectral definition of an expander graph: we want $d - |\lambda_2|$ to remain uniformly bounded away from 0.
These results suggest that $d-|\lambda_2|$ may be a good quantitative way of measuring the quality of an expander graph: the larger, the better. This is complemented by the Alon-Boppana bound \cite{nilli} stating that for any $d-$regular graph on $n$ vertices
$$ |\lambda_2| \geq 2\sqrt{d-1} + o_{n}(1).$$
Here $o_n(1)$ tends to 0 as $n \rightarrow \infty$. This is the sharp bound. Graphs for which
$$ \max(|\lambda_2|, |\lambda_n|)  \leq 2\sqrt{d-1}$$
are called `Ramanujan' and are best possible in terms of spectral expansion.

\subsection{Friedman's result.} A special case of a result of Friedman implies that if the reals $x_1, \dots, x_n$ are i.i.d. samples from the same distribution, then, with high likelihood, the arising graph is close to Ramanujan. 
\begin{thm}[Friedman, Theorem 1.2. in \cite{fried}] If $x_1, \dots, x_n$ are i.i.d. random variables chosen from an absolutely continuous distribution and $G$ is the 4-regular graph constructed from them as above, then for any $\varepsilon> 0$ there exists $c>0$ such that with likelihood at least $1-c/n$, we have
$$ |\lambda| \leq 2 \sqrt{3} + \varepsilon.$$ 
\end{thm}

The formulation in Friedman \cite{fried} is slightly different (because the underlying result is more general): there, $2$ cyclic permutations are used to generate edges on $\left\{1,2,\dots,n\right\}$. Here, we are only dealing with one deterministic cycle (based on their index) and one random cyclic permutation (based on their order) which is easily seen to be equivalent to the
case of two random cyclic permutations as follows: for any two independent random cyclic permutations $\pi_1, \pi_2$, we can use the first to define the cycle and the second as then establishing random connections: for any fixed $\pi_1$ and random $\pi_2$, the graph arising from $(\pi_1, \pi_2)$ behaves exactly as our model. Alternatively, $\pi_1$ serves as a relabeling of the vertices which has no impact on spectral properties.

\section{Testing for Pseudo-Randomness}
Suppose someone hands us a sequence $x_1, \dots, x_n$ and claims that these are i.i.d. samples from a random variable. How would one go about checking such a claim? This leads into the realm of randomness tests of which there are many. A first systematic group of tests was proposed by Knuth \cite{knuth}.  Marsaglia \cite{marscd} produced a group of such tests which became known as the \textsc{diehard tests}. L'Ecuyer \& Simard \cite{lec} produced \textsc{TestU01}, a software library comprised of a large number of tests. 
The point of this paper is to propose a new test.
\begin{quote}
\textbf{The Expander Test}\textbf{.} Generate the graph $G$ and compute $\lambda = \max\left\{ |\lambda_2(A)|, |\lambda_n(A)|\right\}$, where $A$ denotes the adjacency matrix. The larger $\lambda - 2\sqrt{3}$ is, the less likely the sequence is to be random.
\end{quote}

The formulation is intentionally vague: to the best of our knowledge, not much about the distribution of $\lambda - 2 \sqrt{3}$ (which does depend on $n$) is known. In practice, however, it is easy to get an estimate of what the typical for random numbers of a certain cardinality since the quantity concentrates quite quickly: we discuss this in the next section. As usual, passing the test says little about the true nature of the sequence while failing the test is a strong indicator that the sequence is not comprised of i.i.d. samples.

\subsection{Benchmarks.} 
We start with a simple test of the distribution of the second largest eigenvalue for purely random numbers where our source is the random number generator built into \textsc{Mathematica} (which, incidentally, passes the Expander Test). We observe (see Fig. 5) that both distributions, those for $n=1000$ and those for $n=10000$ look quite similar (the second one being more tightly concentrated). Moreover, the mean does seem to approach $2\sqrt{3}$ from below. 
\begin{center}
\begin{figure}[h!]
\begin{tikzpicture}
\node at (0,0) {\includegraphics[width=0.4\textwidth]{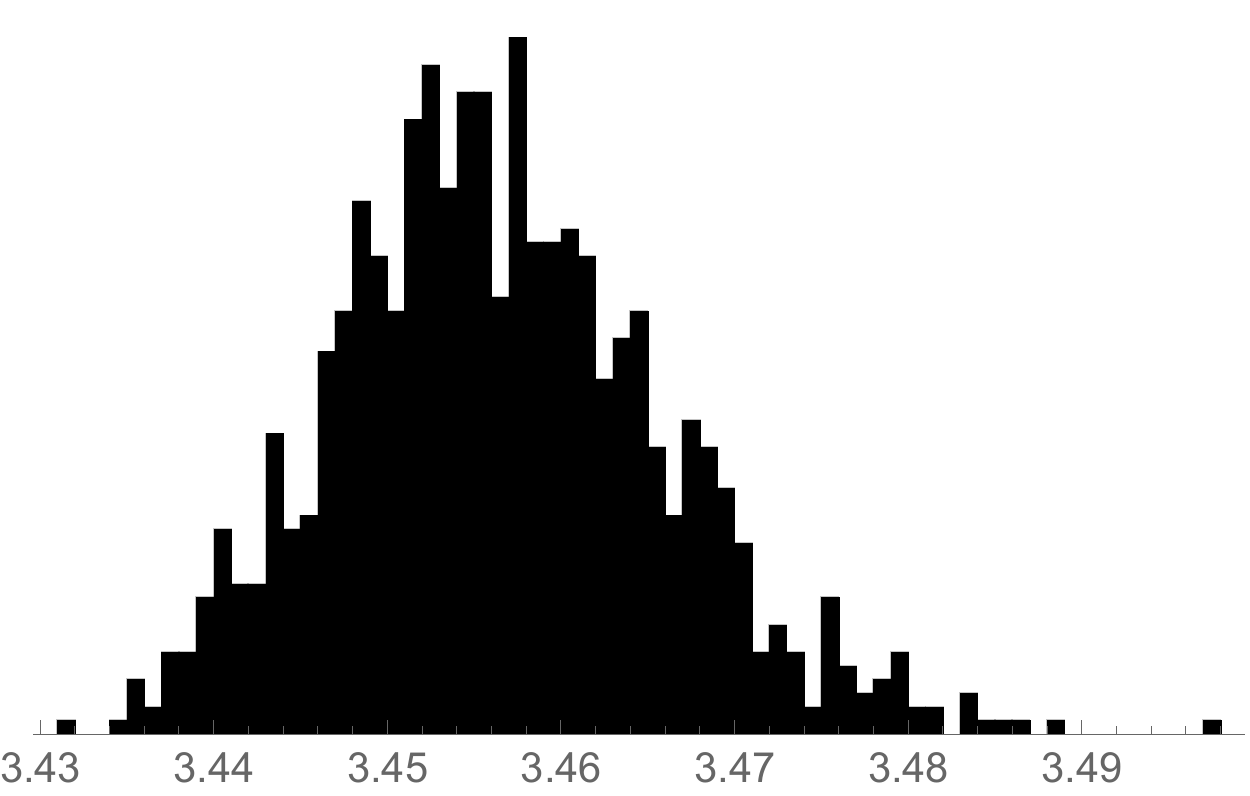}};
\node at (6,0) {\includegraphics[width=0.4\textwidth]{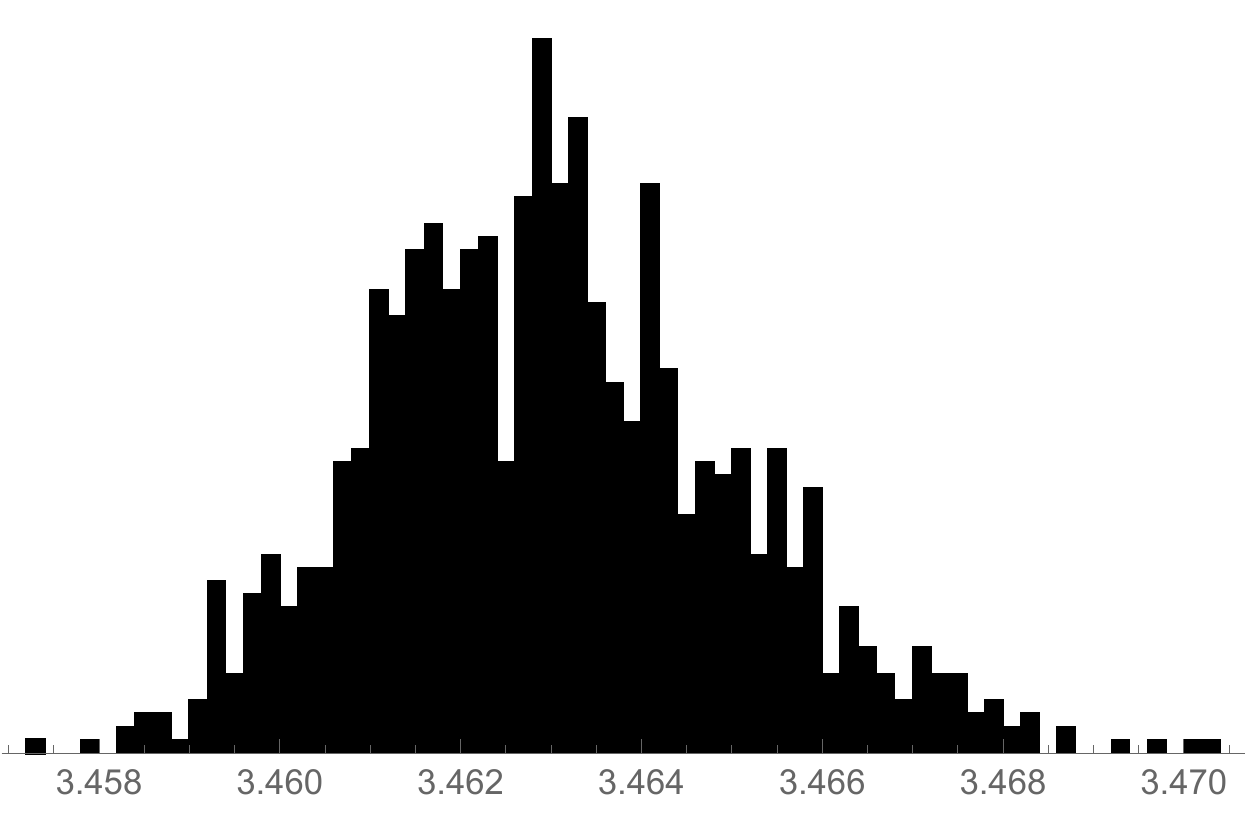}};
\end{tikzpicture}
\caption{Distribution of the second largest eigenvalue generated from $n=1000$ (left) and $n=10000$ (right) random reals.}
\end{figure}
\end{center}

This is also substantiated by our Monte Carlo estimates in Table 1. While more precise theoretical estimates would be quite desirable, these are likely to be difficult (as similar questions for the $d-$regular graphs are still open).

\begin{center}
\begin{table}[h!]
\begin{tabular}{ l | c | c | c | c | c |c}
  $n$ & 10 & 50 &100 & 500 & 1000 & 2000 \\
  \hline
  Mean of $\lambda$ &    3 & 3.36 & 3.41 & 3.45   & 3.457 & 3.46 \\
  Standarddeviation $\lambda$ & 0.295  &0.076 &  0.046  & 0.015 & 0.01 & 0.005  \\
  $\mathbb{P}(\lambda < 2 \sqrt{3})$ & 93\% & 89.5\% & 87.5\% & 81.2\% & 79\% & 78.6\% \\
  \vspace{1pt}
\end{tabular}
\caption{Monte Carlo estimates for properties of Graphs generated from $n$ randomly generated numbers.}
\end{table}
\end{center}

\subsection{Why is this a reasonable test?} We quickly discuss why this is a reasonable test: this is a confluence of several factors
\begin{enumerate}
\item the order in which i.i.d. random variables show up is completely unconnected to their (relative) size and
\item in particular, if we were to partition the vertices of our graph into $S \cup T = \left\{1,2,\dots, n\right\}$, then the
number of edges between $S$ and $T$ should only depend on the cardinality of the sets $|S|, |T|$.
\item Since a search over all possible subsets is not computationally feasible, the largest eigenvector of the adjacency
matrix $A$ can be seen as a reasonable relaxation of this test and
\item thanks to the Alon-Boppana bound  \cite{nilli} and Friedman's theorem \cite{fried}, we have a good theoretical understanding of what to expect: purely random variables will be close to $\lambda \sim 2\sqrt{3}$ while pseudo-random numbers might be close; however, anything with $\lambda \gg 2\sqrt{3}$ away will not be random.
\end{enumerate}

In particular, regarding the eigenvalue as a suitable relaxation, we also mention the Expander Mixing Lemma \cite{ac86}: let $G$ be a $d-$regular graph, let $\lambda = \max(|\lambda_2|, |\lambda_n|)$ and let $S,T \subseteq V$. Then, denoting the number of edges between $S$ and $T$ by $|e(S,T)|$, we have
$$ \left| |e(S,T)| - \frac{d}{n} |S| \cdot |T| \right| \leq \lambda \sqrt{|S| \cdot |T|}.$$
This statement simply states that the expected number of edges between any subsets is close to its expectation. There also exists a converse statement due to Bilu \& Linial \cite{bilu}.  We are not aware of any work in this direction. However, there is work following the opposite direction of reasoning: using expander graphs to generate pseudo-random sequences. 
We refer to Impagliazzo \& Zuckerman \cite{imp} and Lauter, Charles \& Goren \cite{lauter}.

\subsection{Implementation.} A priori, it seems like the test is quite costly to run since, when testing the first $n$ elements, we are interested in eigenvalues of a $n \times n$ matrix. However, there are a couple of helpful facts that simplify life:
\begin{enumerate}
\item the largest eigenvalue is 4 and we know the largest eigenvector (given by the constant vector),  we are only interested in the second largest eigenvalue
\item we only require a \textit{lower} bound on the second eigenvalue to rule out randomness of a sequence
\item and the matrix is both sparse and symmetric.
\end{enumerate}
This makes power iteration methods a very reasonable tool: start with a random vector $x_0$ whose mean value is 0 and define
$ x_{k+1} = Ax_k/\|A x_k\|.$
Since the vector $x_k$ is orthogonal to the constant vector (corresponding to the largest eigenvector), so is $A x_k$. (In practice, it is a good idea to occasionally normalize again to compensate for numerical errors). Since $A$ is sparse, matrix multiplication is very cheap. $| \left\langle x_k, A x_k\right\rangle$ then serves as lower bound on the second largest eigenvalue. With this method, even matrices of size $n=50.000$ do not pose a serious challenge.

\subsection{A Concluding Question.} It would be quite interesting to understand whether the Expander Test can be extended: are there other ways of naturally associating graphs to sequences of numbers such that randomness in the sequence is being reflected in reasonable expansion properties?

\section{Examples}
Throughout this section, we will use the notation $\lambda(G_n)$ to denote the second largest eigenvalue of the Graph computed using the first $n$ elements of the sequence. Since all our examples here are integer-valued and our construction requires the numbers to be distinct, we add a tiny random variable to each element in the sequence (this corresponds to fixing a random order on the occurence of random elements) -- since this can only increase randomness, it is a valid procedure.

\subsection{Lehmer's original example.} Lehmer's \cite{lehmer} sequence is considered the earliest (1949) published method for generating pseudorandom numbers. It is given by $x_1 = 47594118$ ('chosen at random from a wastepaper basket of punched cards')
\begin{align*}
x_{n+1} = 23 \cdot x_n \mod 10^{8}+1.
\end{align*}
The sequence is considered to be somewhat reasonable Knuth \cite{knuth} records it as passing $\chi-$square tests without difficulty. However, it is also understood that the multiplier 23 is somewhat small. Running our statistical tests on the graph arising from the first 500 elements, we already obtain $\lambda(G_{500}) = 3.54$. This gets worse as the number of elements increases, for example $\lambda(G_{2000}) = 3.55$ and $\lambda(G_{5000}) = 3.57$. These numbers strongly indicate that the sequence is not truly random.\\

\subsection{Examples suggested by Park-Miller \cite{park}.} Lehmer's example is representative of a larger group of random number generators of the type
$$ x_{n+1} = a\cdot x_{n} + c \mod m.$$
For these examples, one has to be careful in choosing appropriate $a,c,m$, this is discussed at great length in Knuth \cite{knuth}. Park \& Miller \cite{park} collect a number of examples from the literature where proposed parameters do \textit{not} result in suitable sequences (i.e. failing various statistical tests). One such sequence is  
$$ x_{n+1} = 20403 \cdot x_{n} \mod 2^{15}.$$
The period of this sequence is rather short (only $2^{13} \sim 8192$). However, we see problems much earlier since $\lambda(G_{1000}) = 3.471$, $\lambda(G_{2000}) = 3.48$, $\lambda(G_{3000}) = 3.496$ and $\lambda(G_{5000}) = 3.533$. While the first two values are somewhat acceptable, the others are outside what one would consider reasonable. Another example is a sequence proposed in a variety of textbooks  `and is now
something of an emerging standard in the undergraduate computer science textbook market' \cite{park}. It reads
$$ x_{n+1}= 25173 \cdot x_{n} + 13840 \mod 2^{16}.$$
We see that $\lambda(G_{5000}) = 3.58$, easily dismissing the sequence as not random. We continue with an example that was implemented as the random number generator in the programming language \textsc{Turbo Pascal}
$$ x_{n+1} = 129 \cdot x_{n} + 907633385 \mod 2^{32}.$$
We observe $\lambda(G_{5000}) = 3.508$ and $\lambda(G_{10000}) = 3.52$ which easily dismisses the sequence. Rotenberg \cite{rot} suggested to use
$$ x_{n+1} = (2^7+1) \cdot x_{n} + 1 \mod 2^{35}.$$
Knuth \cite{knuth} refers to this sequence as `borderline' with regards to its behavior under the $\chi^2-$test. It fails our test dramatically since $\lambda(G_{5000}) = 3.54$.
All the examples above tend to emulate uniform distribution; this is not necessary for our test -- we care only about whether they come from some distribution in an i.i.d. fashion. A simple example is the logistic map, we pick $x_1 = 0.3$ and
$x_{n+1} = 3.98 x_n (1-x_n).$
We observe $\lambda(G_{100}) = 3.75$, this is clearly (and unsurprisingly) not an i.i.d. sequence.  

\subsection{Some `good' examples.}
We now discuss some sequences that perform quite well. Knuth \cite{knuth} discusses the `satisfactory' example (with regards to the $\chi^2-$test)
\begin{align*}
x_1 &= 0\\
x_{n+1} &= 3141592653 \cdot x_n + 2718281829 \mod 2^{35}
\end{align*}
We compute $\lambda(G_{10000}) = 3.465$ and $\lambda(G_{20000}) = 3.463$ showing that it indeed behaves quite well with regards to our test. Lewis, Goodman \& Miller \cite{lewis} propose
$$x_{n+1} = x_n \cdot 16807 \mod 2147483647$$
which Park \& Miller \cite{park} describe as a good minimal standard. We compute
$\lambda(G_{20000}) = 3.463$. L'Ecuyer \cite{lec} proposes
$$ x_{n+1} = x_n \cdot 40692 \mod (2^{31} -249),$$
we observe $\lambda(G_{20000}) = 3.461$.
Knuth \cite{knuth} mentions a method proposed by R. Coveyou which is known to have long periods
\begin{align*}
x_1 &\mod 4 =2\\
x_{n+1} &= x_n(x_n+1) \mod 2^{e}
\end{align*}
We run this method with $x_1=1234$ and $e=32$. We find that $\lambda(G_{10000}) = 3.462$. We conclude with the lagged Fibonacci sequence (see e.g. Knuth \cite{knuth}) which is initialized with 56 integers and then given by
$$ x_{n} = x_{n-24} + x_{n-55} \mod 2^k.$$
It is known to have period $2^{k-1}(2^{55}-1)$. We found the sequence to uniformly yield good values; this agrees with Marsaglia \cite{marspeech} who noted that the sequence passes a large number of tests.\\

\textbf{Acknowledgments.}
The author is grateful to Noah Kravitz for helpful conversations.


\begin{thebibliography}{10}

\bibitem{alo86} N. Alon. Eigenvalues and expanders. Combinatorica, 6(2):83--96, 1986.

\bibitem{ac86} N. Alon and F. R. K. Chung. Explicit construction of linear sized tolerant networks,
Discrete Math., 72:15--19, 1989

\bibitem{am85} N. Alon and V. D. Milman. $\lambda_1$, isoperimetric inequalities for graphs, and superconcentrators. J. Combin. Theory Ser. B, 38(1985), p. 73--88

\bibitem{bo} J. Bourgain and P. Varju, Expansion in SLd(Z/qZ), q arbitrary, Inventiones Mathematicae 188 (2012): p. 151--173

\bibitem{bgt3} E. Breuillard, B. Green and T. Tao, Suzuki groups as expanders, Groups Geom. Dyn.
5 (2011), 281--299.

\bibitem{bilu} Y. Bilu and N. Linial. Lifts, discrepancy and nearly optimal spectral gaps. Combinatorica 26 (5) (2006) p. 495--519

\bibitem{bus} P. Buser. A note on the isoperimetric constant. Ann. Sci. ENS 
15 (1982), p. 213--230. 

\bibitem{che} J. Cheeger. A lower bound for the smallest eigenvalue of the Laplacian. In Problems
in analysis (Papers dedicated to Salomon Bochner, 1969), pages 195--199. Princeton
Univ. Press, Princeton, NJ, 1970

\bibitem{dod84} J. Dodziuk. Difference equations, isoperimetric inequality and transience of certain random walks. Trans. Amer. Math. Soc., 284 (1984):787--794.

\bibitem{lec} Pierre L’Ecuyer and Richard Simard, TestU01: A Software Library in ANSI C for Empirical Testing of Random Number Generators, ACM Trans. Mathematical Software 33 (2007): 22.

\bibitem{eich} Jurgen Eichenauer and Jurgen Lehn, A non-linear congruential pseudo random number
generator, Statistische Hefte 27, 315-326 (1986)

\bibitem{fried} J. Friedman, A proof of Alon’s second eigenvalue conjecture and related problems,
Mem. Amer. Math. Soc. 195 (2008), no. 910

\bibitem{linial} S. Hoory, N. Linial and A. Widgerson, Expander Graphs and their Applications, Bull. Amer. Math. Soc. 43 (2006), p. 439--561.

\bibitem{imp} R. Impagliazzo and D. Zuckerman, How to recycle random bits, 30th Annual Symposium on Foundations of Computer Science, 1989

\bibitem{k2}  M. Kassabov, Symmetric groups and expander graphs, Inventiones Mathematicae 170 (2007), no.
2, 327--354.

\bibitem{knuth} D.E. Knuth, The Art of Computer Programming: Seminumerical Algorithms, Vol. 2, 2nd ed.
(Addison-Wesley, 1981).

\bibitem{kol} A. N. Kolmogorov and Y.M. Barzdin, On the realization of nets in 3-dimensional
space, Probl. Cybernet, 8, 261-268, 1967. See also Selected Works of A.N. Kolmogorov,
Vol. 3, pp. 194-202 (and a remark on page 245), Kluwer Academic Publishers, 1993

\bibitem{kow} E. Kowalski, An Introduction to Expander Graphs,  Cours Specialises Volume: 26, Societe Mathematique de France, 2019.

\bibitem{lauter} K. Lauter, D. Charles and E. Goren, Pseudorandom number generation with expander graphs, United States Patent, US7907726B2, 2006


\bibitem{lec} P. L'Ecuyer, Efficient and portable combined random number generators. Communications of the ACM 31 (1988), p. 742--749.


\bibitem{lehmer} D. H. Lehmer: Mathematical methods in large-scale computing units, in Proceedings of the 2nd Symposium on Large-Scale Digital Calculating Machinery, Cambridge, MA, 1949, pp. 141--146, Cambridge, MA, 1951, Harvard University Press.


\bibitem{lewis} P. Lewis; A. Goodman, J. Miller, A pseudo-random number generator for the system/360. IBM Systems Journal. 8 ( 1969):  p. 136--143.


\bibitem{lub} A. Lubotzky, Expander Graphs in Pure and Applied Mathematics, Bull. Amer. Math. Soc. 49 (2012), p. 113 -- 162.

\bibitem{lps1}  A. Lubotzky, R. Phillips and P. Sarnak, Ramanujan conjecture and explicit construction of expanders, Proc. STOC. 86 (1986), 240--246.

\bibitem{lsv2}  A. Lubotzky, B. Samuels and U. Vishne, Explicit constructions of Ramanujan complexes of type $\tilde A_d$, European J. Combin. 26 (2005), no. 6, 965--993.

\bibitem{m1} G.A. Margulis, Explicit constructions of expanders. (Russian) Problemy Peredaci Informacii 9 (1973), no. 4, 71--80. English translation: Problems of Information Transmission 9 (1973), no. 4, 325--332 (1975).

\bibitem{m2} G. Margulis, Explicit constructions of graphs without short cycles and low density
codes, Combinatorica 2 (1982), no. 1, 71--78.

\bibitem{m3} G, Margulis, G.A. Margulis, Explicit group-theoretic constructions of combinatorial schemes and
their applications in the construction of expanders and concentrators, Problems of
Information Transmission, 24(1):39--46, 1988.

\bibitem{mars} G. Marsaglia, Random numbers fall mainly in the planes. PNAS 61 (1968), p. 25--28. 

\bibitem{marspeech} G. Marsaglia,  A Current View of Random Number Generators, Proceedings of `Computer Science and Statistics: 16th Symposium on the Interface'  Atlanta, Elsevier, 1984.

\bibitem{marscd} G. Marsaglia, The Marsaglia Random Number CDROM including the Diehard Battery of Tests of Randomness. Florida State University. 1995.

\bibitem{modula}  MacModula-2 System Reference Manual. Modula Corporation. 1985. p. 41.


\bibitem{mo} M. Morgenstern, Existence and explicit constructions of q + 1 regular Ramanujan graphs for every prime power q, J. Combin. Theory Ser. B 62 (1994), 44--62

\bibitem{nilli}  A. Nilli. On the second eigenvalue of a graph. Discrete Math., 91(2):207--210, 1991.

\bibitem{park} S. Park and K. Miller, Random Number Generators: Good Ones Are Hard to Find, Communications of the ACM 31 (1988): p. 1192--1201.

\bibitem{pinsker} M. S. Pinsker. On the complexity of a concentrator. In 7th International Telegraffic
Conference, pages 318/1--318/4, 1973.

\bibitem{ps2} L. Pyber and E. Szabo, Growth in finite simple groups of Lie type of bounded rank,  J. Amer. Math. Soc. 29 (2016), 95--146

\bibitem{rvw} O. Reingold, S. Vadhan and A. Wigderson, Entropy waves, the zig-zag graph product,
and new constant-degree expanders, Ann. of Math. (2) 155 (2002), no. 1, 157--187.

\bibitem{rot} A. Rotenberg, A new pseudo-random number generator. Journal of the ACM (JACM) 7.1 (1960): 75--77.

\bibitem{va2}  A. Valette, Graphes de Ramanujan et applications, (French) Ramanujan graphs and
applications, Seminaire Bourbaki, Vol. 1996/97. Ast´erisque No. 245 (1997), Exp. No.
829, 4, p. 247--276


\end{thebibliography}
\end{document}